\documentclass{article}

\usepackage[utf8]{inputenc}
\usepackage[T1]{fontenc}
\usepackage{lmodern}

\usepackage{amsmath, amssymb, amsfonts}

\usepackage[a4paper, margin=3cm]{geometry}

\usepackage{graphicx}
\usepackage{subcaption}
\usepackage{float}

\usepackage{xcolor}
\usepackage{tcolorbox}

\usepackage{authblk}

\usepackage[hidelinks,unicode=true]{hyperref}
\usepackage{parskip}

\newtheorem{lemma}{Lemma}[section]
\newtheorem{Proposition}[lemma]{Proposition}

\newtheorem{proof}{Proof}

\title{Importance sampling for Sobol' indices estimation}

\author[1,2]{Haythem Boucharif}
\author[1,2]{Jérôme Morio}
\author[2]{Paul Rochet}

\affil[1]{ONERA/DTIS, Université de Toulouse, Toulouse, 31055, France}
\affil[2]{Fédération ENAC ISAE-SUPAERO ONERA, Université de Toulouse, Toulouse, France}

\date{March 01, 2026}

\begin{document}

\maketitle

\begin{abstract}
We propose a new importance sampling framework for the estimation and analysis of Sobol' indices.
We focus on the estimation of the conditional second-moment quantity underlying these indices, which is the most challenging term to estimate.
We show that this quantity, originally defined under a reference input distribution, can be estimated from samples drawn under auxiliary distributions by reweighting the model outputs.
We derive the optimal sampling distribution that minimises the asymptotic variance of efficient estimators and demonstrate its impact on estimation.
Beyond variance reduction, the framework also supports distributional sensitivity analysis through reverse importance sampling.
\end{abstract}


\section{Introduction}
Complex computer models are now widely used across many scientific and engineering fields. Since they are often costly to evaluate, handling uncertainty and performing sensitivity analysis is particularly critical. In the context of global sensitivity analysis (GSA), where model inputs $\mathbf{X} = (X_1, \dots, X_k)$ are treated as random variables, the aim is to understand, with as few model evaluations as possible, the global influence of each input or group of inputs on the variability of the model output $Y$ \cite{Saltelli2000,Sobol1993}.

To quantify this impact, a large variety of sensitivity indices have been proposed in the literature, each one associated with a different criterion. Some indices rely on higher-order moments \cite{Owen2013,OwenDickChen2014}, others are based on discrepancy measures between probability distributions \cite{Borgonovo2007,Borgonovo2011,DaVeiga2015}, while others use contrast-based functionals or Cramér--von Mises distance \cite{Fort2016,Gamboa2018}. These indices allow for a refined view of the output distribution. Among all these tools, Sobol' indices remain the most popular in practice. Originally introduced in \cite{Sobol1993}, they are particularly well suited for scalar outputs, and can be generalised to vector-valued or functional outputs. Thanks to the Hoeffding decomposition \cite{Hoeffding1948}, Sobol' indices measure the contribution of each group of inputs to the output variance. The first-order Sobol' indices quantify the main effect of each input, while the total-order Sobol' indices also account for interaction effects. Hence, they are variance-based measures that capture the contribution of inputs to the overall variability of the model output.\\

In most practical applications, the model is a complex black-box code whose output distribution is unknown. Analytical formulas for Sobol’ indices are thus rarely available. As a result, these indices must be estimated, making their statistical estimation a key challenge in sensitivity analysis, especially for the conditional second moment $\eta = \mathbb{E}[\mathbb{E}[Y \mid X]^2]$. Over the last decades, many estimation strategies have been proposed and can be broadly classified into four categories: (i) Monte Carlo and quasi-Monte Carlo methods, including plain random sampling, low-discrepancy sequences, and nested designs \cite{Goda2017,Owen2013Small}; (ii) spectral and expansion techniques, such as the Fourier Amplitude Sensitivity Test (FAST) \cite{Cukier1978}, Random Balance Designs (RBD) \cite{Tarantola2006}, the Effective Algorithm for Sensitivity Indices (EASI) \cite{Plischke2010}, and polynomial chaos expansions \cite{Sudret2008}; (iii) the Pick–Freeze procedure, which directly estimates conditional variances via sample-based manipulations \cite{Gamboa2016,Janon2014}; and (iv) nonparametric estimators including nearest-neighbor \cite{Devroye2003,Devroye2018,Broto2020} 
and kernel-based methods \cite{DaVeiga2013,Gamboa2020,DaVeiga2024}, 
and the rank-based estimator \cite{KleinLagnouxRochetNguyen2024, Gamboa2022} 
built on Chatterjee’s correlation coefficient \cite{Chatterjee2021} which share the practical advantage of being applicable in the so-called given-data setting, where only a fixed input–output sample set is available and direct access to the black-box model is not possible.

The theoretical properties of methods from categories (iii) and (iv) are now well established. In particular, for first-order Sobol’ indices, consistency and asymptotic normality have been proven for kernel estimators \cite{DaVeiga2024,DaVeiga2013}, the Pick–Freeze estimator \cite{Janon2014}, nearest-neighbor methods \cite{Devroye2018}, and rank-based estimators \cite{Gamboa2022}.\\

With the increasing complexity of models and growing computational costs, the Pick--Freeze setting is becoming less practical in many applications, since the number of model evaluations required to maintain a given accuracy for all indices increases linearly with the number \( k \) of input variables. When the context allows it, using this particular design can considerably simplify estimation, but the number of code calls often makes it prohibitive. By contrast, given-data methods are increasingly favored because they operate from a provided input--output sample only. This naturally raises the question of efficiency limits in the given-data setting.

Among them, several estimators achieve the asymptotic efficiency lower bound. 
The expression of the efficiency bound for $\eta$ was first derived in \cite{DoksumSamarov1995} for a continuous input $X$, 
from the efficient influence function associated with the parameter. 
The same work also proposed an efficient estimator of a truncated version of $\eta$. 
Since then, several methods have been shown to be asymptotically efficient in the given-data setting. 
For first-order indices, the kernel estimator introduced in \cite{DaVeiga2013} achieves the asymptotic optimal variance bound for estimating $\eta$ (see Theorems~3.4 and~3.5). 
This result was initially established for first-order Sobol' indices and later extended to higher-order ones under suitable regularity conditions \cite{DaVeiga2024} (see Theorems~4.1 and~4.2). 
These estimators, however, remain difficult to implement in practice. 
In contrast, the rank-based estimator of \cite{Gamboa2022}, built on order statistics, is much simpler but sub-optimal and restricted to first-order indices, 
although it has been further generalised to attain the efficiency bound in \cite{KleinRochet2024} (Proposition 3.3). 
The question of efficiency bounds has been thoroughly discussed in \cite{KleinLagnouxRochetNguyen2024}. 
In this work, we mainly focus on these efficient estimators.\\

Even though these estimators are already asymptotically optimal under the reference input distribution, their accuracy may degrade when the sample size is small or when the model response is highly nonlinear.
In such cases, a large fraction of the simulation effort may be wasted on regions of the input space that contribute little to the Sobol’ quantity of interest, while the most informative regions are poorly explored. When the black-box model is still available, a natural way to improve efficiency is to modify the sampling scheme so that evaluations are more concentrated on informative regions.
This motivates the use of \emph{importance sampling} (IS) \cite{Kahn1951}, one of the most classical and effective approaches: it modifies the sampling distribution to obtain a more accurate estimator of the target quantity. The variance achieved by IS strongly depends on how this sampling distribution is chosen.\\

In this work, we reexamine the estimation of Sobol' indices from the perspective of importance sampling. More specifically, we focus on the quantity $\eta$, which appears as a key component in the variance of Sobol’ estimators \and is the most challenging term to estimate accurately. We introduce a novel reformulation of $\eta$ using an auxiliary probability distribution,which makes it possible to characterise optimal sampling distribution for efficient estimators, potentially achieving zero variance under appropriate conditions.\\

To the best of our knowledge, no prior work has established a theoretical connection between variance reduction through importance sampling and the estimation of Sobol' indices. Two related contributions deserve to be mentioned. First, \cite{Beaurepaire2014} propose a computational approach sharing similarities with importance sampling: they introduce an auxiliary variable correlated with a given input and estimate conditional expectations using a Gauss–Hermite quadrature scheme. Their main contribution lies in the simultaneous estimation of all first-order Sobol’ indices using only a single Monte Carlo simulation, thereby reducing numerical effort compared to classical procedures. However, their method does not address variance minimisation explicitly. Second, \cite{Borgonovo2024} introduce a weighting method for estimating Sobol’ total effects under constrained input spaces, where an importance-type density ratio is used to correct Jansen’s estimator in non-Cartesian domains. In contrast, our contribution focuses on the classical aspect of importance sampling for variance reduction and provides a fully flexible probabilistic framework that also extends naturally to distributional sensitivity analysis.

Beyond variance reduction, the proposed importance sampling framework can also be exploited to assess the robustness of Sobol' indices with respect to changes in the input distributions. 
This corresponds to the concept of \emph{reverse importance sampling} \cite{Beckman1987,BadeaBolardo2008,Lemaitre2015}, which exploits the fact that a single dataset can be reused to evaluate sensitivity measures under alternative probabilistic scenarios. 
Note that this should not be confused with the reverse importance sampling estimator (or harmonic mean estimator) used in Bayesian model evidence estimation \cite{NewtonRaftery1994}. A post-processing importance sampling approach was also proposed in \cite{Sparkman2016}, where existing Monte Carlo samples are reweighted to estimate global sensitivity indices. Their method relies on a kernel-like weighting scheme to approximate conditional moments.\\

The paper is structured as follows. Section 2 introduces Sobol’ indices and the classical importance sampling framework. Section 3 presents our main theoretical results. Section 4 illustrates these results with numerical experiments, including variance reduction analysis and a reverse importance sampling study of distributional sensitivity on a real-life test case. Section 5 concludes and outlines future research directions.

\section{Main notation}
\subsubsection*{Definition of Sobol' indices}
Consider a triplet $(Y, \mathbf{X}, \mathbf{W})$  of random variables satisfying the relation $Y = f(\mathbf{X},\mathbf{W})$ for some measurable function \(f\), where \(Y\) is a real-valued square-integrable output, \(\mathbf{X} = (X_1, \dots, X_k) \in \mathbb{R}^k\) is a vector-valued input, and \(\mathbf{W}\) is a random vector independent of \(\mathbf{X}\). In the context of sensitivity analysis, the so-called Sobol' index $S_u$ is used to quantify the impact of a subset \(u \subseteq \{1, \dots, k\}\) of the input variables on the output \(Y\), and is defined as:
\begin{equation}
S_u = \frac{\mathrm{Var}\!\left(m(\mathbf{X}_u)\right)}{\mathrm{Var}(Y)},
\end{equation}

where \( m(\mathbf{x}_u) = \mathbb{E}[Y \mid \mathbf{X}_u = \mathbf{x}_u] \), with 
\(\mathbf{X}_u = (X_j)_{j \in u}\) the subvector of inputs corresponding to the index set \(u\), 
and \(\mathbf{X}_{-u} = (X_j)_{j \notin u}\) its complementary subvector containing the remaining 
components of \(\mathbf{X}\). One can develop the variance term as 
\(\mathrm{Var}(m(\mathbf{X}_u)) = \mathbb{E}[m^2(\mathbf{X}_u)] - (\mathbb{E}[Y])^2\), and both expectations and variances can in practice be estimated through classical Monte Carlo estimators based on independent samples of \((Y,\mathbf{X})\).
Hence, as recalled above, the main challenge for inference purposes lies in handling the conditional expectation \(m(\mathbf{X}_u)\) itself in order to construct an estimator of

\begin{equation}
\eta_u = \mathbb{E}\!\left[m^2(\mathbf{X}_u)\right]
 = \mathbb{E}\!\left[m(\mathbf{X}_u)\,\mathbb{E}[Y\mid \mathbf{X}_u]\right]
 = \mathbb{E}[\,\mathbb{E}[Y\,m(\mathbf{X}_u)\mid \mathbf{X}_u]\,]
 = \mathbb{E}\!\left[Y\, m(\mathbf{X}_u)\right],
\end{equation}
\noindent
where the second equality follows from the law of total expectation.

For an i.i.d. sample $(\mathbf{X}^{(i)}, Y^{(i)}),\; i=1,\dots,n$, 
efficient estimators of $\eta_u$ are asymptotically Gaussian, 
with variance achieving the lower bound \cite{DoksumSamarov1995}. Writing $\nu(\mathbf{X}_u)=\mathrm{Var}(Y\mid\mathbf{X}_u)$,
the optimal asymptotic variance is given by:
\begin{equation}
\sigma^2_{\mathrm{opt},u}
=4\,\mathbb{E}\!\left[m^2(\mathbf{X}_u)\,\nu(\mathbf{X}_u)\right]
+\mathrm{Var}\!\left(m^2(\mathbf{X}_u)\right).
\label{eq:sigma_opt_rewrite}
\end{equation}

As motivated in the introduction, we focus on improving the estimation of the key quantity $\eta_u$. 
For both theoretical and practical reasons, it is often beneficial to modify the original distribution of the input variables. 
Our goal is to investigate whether it is possible to reduce the variance of the estimator of $\eta_u$ by judiciously choosing a sampling distribution $q$.\\

In the importance sampling framework, briefly recalled below, the input vector $\mathbf{X}$ is initially drawn from a reference distribution $p$, under which the components may be dependent, although independence can be assumed for analytical convenience. 
We then introduce an alternative sampling distribution $q$, possibly with a different dependency structure. 
To avoid ambiguity, expectations are indexed by the underlying distribution: for instance, $\mathbb{E}_p[\cdot]$ denotes expectation under $p$. 
The distribution of $\mathbf{W}$ is fixed and remains unchanged, and the output $Y$ is affected only through changes in the distribution of $\mathbf{X}$.

\subsubsection*{Importance sampling} 
This classical variance-reduction method estimates expectations of the form $\mathbb{E}_p[f(\mathbf{X})]$ 
by rewriting them under another distribution $q$ as $\mathbb{E}_q[f(\mathbf{X})\, w(\mathbf{X})]$,
where $w(\mathbf{x}) = p(\mathbf{x}) / q(\mathbf{x})$ denotes the likelihood ratio.
To ensure unbiasedness, the support of $q$ must include that of $f(\mathbf{x})\, p(\mathbf{x})$. 
Given an i.i.d. sample $\mathbf{X}^{(1)}, \dots, \mathbf{X}^{(n)}$ drawn from $q$, the corresponding estimator is
\begin{equation}
    \hat{I}^{\mathrm{IS}} = \frac{1}{n} \sum_{i=1}^n f(\mathbf{X}^{(i)})\, w(\mathbf{X}^{(i)}).
\end{equation}
This estimator achieves zero-variance if and only if $f$ is nonnegative and $q(\mathbf{x}) \propto f(\mathbf{x})\,p(\mathbf{x})$ almost everywhere.
This ideal distribution cannot be used in practice because it depends on the unknown quantity to be estimated. 
However, adaptive procedures exist to approximate it, such as nonparametric methods \cite{zhang1996} or parametric approaches like the cross-entropy method \cite{papaioannou2019improved, deboer2005tutorial}. 
The efficiency of IS \cite{OwenZhou2000} critically depends on the overlap between the support of $q$ and the regions where $p(\mathbf{x}) f(\mathbf{x})$ is large. 
Poorly chosen auxiliary distributions can lead to high-variance or even unstable estimators due to extreme weight values.\\

The next section extends the classical IS framework to characterise the optimal $q$ for estimating $\eta_u$ and improving the estimation of Sobol’ indices.

\section{Main results}
\label{sec:main}
\subsection{Input distribution change via importance sampling} 

We introduce the importance sampling framework by allowing a modification of the input distribution of the random vector $\mathbf{X} \in \mathbb{R}^k$, 
while keeping the noise variable $\mathbf{W}$ fixed.
Let $p(\mathbf{x})$ denote the reference joint density of $\mathbf{X}$. 
For any subset of indices $u \subseteq \{1, \dots, k\}$, we let $\mathbf{X}_u = (X_j)_{j \in u}$ and denote by $p_u(\mathbf{x}_u)$ its marginal density. 
We introduce an auxiliary joint density $q(\mathbf{x})$, factorised as
\begin{equation}
q(\mathbf{x}) = q_u(\mathbf{x}_u)\, q_{-u \mid u}(\mathbf{x}_{-u} \mid \mathbf{x}_u),
\end{equation}
where $q_u(\mathbf{x}_u)$ and $q_{-u \mid u}(\cdot \mid \mathbf{x}_u)$ denote the marginal and conditional components of $q$. The corresponding importance weights are then given by
\begin{equation}
w(\mathbf{x}) = \frac{p(\mathbf{x})}{q(\mathbf{x})}
= \frac{p_u(\mathbf{x}_u)\, p_{-u \mid u}(\mathbf{x}_{-u}\mid \mathbf{x}_u)}{q_u(\mathbf{x}_u)\, q_{-u \mid u}(\mathbf{x}_{-u} \mid \mathbf{x}_u)}
= w_u(\mathbf{x}_u)\, w_{-u \mid u}(\mathbf{x}_{-u} \mid \mathbf{x}_u).
\end{equation}

\begin{lemma}\label{lem:Z_u_property}
\textit{
Let $q = q_u\, q_{-u \mid u}$ be any importance sampling distribution on the input $\mathbf{X}$, and define
\begin{equation}
Z_u := \sqrt{w_u(\mathbf{X}_u)}\, w_{-u \mid u}(\mathbf{X}_{-u} \mid \mathbf{X}_u)\, Y.
\end{equation}
Then:
\begin{equation}
\mathbb{E}_q[Z_u \mid \mathbf{X}_u] = \sqrt{w_u(\mathbf{X}_u)}\, m(\mathbf{X}_u).
\end{equation}
In particular, the following identity holds:
\begin{equation}
\eta_u = \mathbb{E}_p\!\left[(\mathbb{E}_p[Y \mid \mathbf{X}_u])^2\right]
= \mathbb{E}_q\!\left[(\mathbb{E}_q[Z_u \mid \mathbf{X}_u])^2\right]
= \mathbb{E}_q\!\left[Z_u\, \mathbb{E}_q[Z_u \mid \mathbf{X}_u]\right].
\end{equation}}
\end{lemma}

The proof is given in the Appendix. \\

In the following, this identity is used to express $\eta_u$ under $q$ and to show that an appropriate choice of $q$ can reduce the variance of its estimator.

\subsection{Optimal asymptotic variance expression using importance sampling}

We define $\phi(\mathbf{x}) = \sqrt{\mathbb{E}[Y^2 \mid \mathbf{X} = \mathbf{x}]}$, noting that $\phi$ is model-specific and does not depend on the choice of sampling distribution.
Under $p$, the optimal asymptotic variance for estimating $\eta_u$ is given by
\begin{equation}
\sigma^2_{\mathrm{opt},u}(p)
=\mathbb{E}_{p}\Big[
4\,m^2(\mathbf{X}_u)\,\nu_p(\mathbf{X}_u)
+m^4(\mathbf{X}_u)
\Big]
-\eta_u^2,
\label{eq:sigma_opt_p_simple}
\end{equation}
where the conditional variance $\nu_p$ is given by
\[
\nu_p(\mathbf{X}_u)
:=\mathrm{Var}_p(Y\mid \mathbf{X}_u)
=\mathbb{E}_p[Y^2\mid \mathbf{X}_u]-m^2(\mathbf{X}_u)
=\mathbb{E}_p[\phi^2(\mathbf{X})\mid \mathbf{X}_u]-m^2(\mathbf{X}_u),
\]

We next extend this result to an auxiliary distribution $q$ and compare it with the reference case.

\begin{Proposition}
\textit{
Let $Z_u$ be defined as in Lemma~\ref{lem:Z_u_property}. 
Then, for the sampling distribution $q$, 
the optimal asymptotic variance is given by}
\begin{equation}
\sigma^2_{\mathrm{opt},u}(q)
=\mathbb{E}_{q}\Big[
4\,\big(\mathbb{E}_{q}[Z_u\mid\mathbf{X}_u]\big)^2
\,\mathrm{Var}_{q}(Z_u\mid\mathbf{X}_u)
+\big(\mathbb{E}_{q}[Z_u\mid\mathbf{X}_u]\big)^4
\Big]
-\eta_u^2.
\end{equation}

\textit{Substituting the expressions of $Z_u$, $\sigma^2_{\mathrm{opt}, u}(q)$ can be equivalently expressed under $p$ as}
\begin{equation}
\sigma^2_{\mathrm{opt},u}(q)
=\mathbb{E}_{p}\Big[
w_u(\mathbf{X}_u)
\big(
4\,m^2(\mathbf{X}_u)\,\nu_q(\mathbf{X}_u)
+m^4(\mathbf{X}_u)
\big)
\Big]
-\eta_u^2,
\label{eq:sigma_opt_q_final}
\end{equation}
\textit{where}
\[
\nu_q(\mathbf{X}_u)
=\mathbb{E}_{p}\big[
w_{-u\mid u}(\mathbf{X}_{-u}\mid \mathbf{X}_u)\,
\phi^2(\mathbf{X})
\mid \mathbf{X}_u
\big]
-m^2(\mathbf{X}_u).
\]
\end{Proposition}

This expression shows that the asymptotic variance depends explicitly on $q$ through the weights $w_u$ and $w_{-u\mid u}$, and can be reduced by optimising them, as will be shown in the next section.

\subsection{Optimal distribution: sequential construction}
We aim to minimise  $\sigma^2_{\mathrm{opt}, u}(q)$ by optimising the components of $q$ in a sequential manner.
The procedure involves two steps. First, for a fixed marginal $q_u$, we determine the optimal conditional density $q^*_{-u \mid u}(\mathbf{x}_{-u} \mid \mathbf{x}_u)$  that minimises the variance contribution given $\mathbf{x}_u$, which is then substituted into the variance expression. In the second step, the resulting objective is minimised with respect to $q_u$.

To characterise the optimal distribution $q^*$, 
a principled approach consists in using the \textit{calculus of variations}.
Following the classical Lagrangian minimisation framework, one can formally derive 
the stationary condition satisfied by $q^*$ under the normalisation constraint. 
In general, verifying that this stationary point corresponds to a minimum would require 
checking second-order optimality conditions. 
However, as mentioned by Owen~\cite[Chapter~9, End Notes (``Calculus of variations''), p.~38]{Owen2013book}, 
the Cauchy--Schwarz inequality is sufficient to establish global minimality in this setting. 
For simplicity, we rely on the following key variational lemma, a standard result in importance sampling theory.

\begin{lemma}\label{lem:variational_minimiser}
\textit{
Let \( \mathcal{D}_k \) denote the set of probability density functions on \( \mathbb{R}^k \).
Then, for any \( g \in \mathcal{D}_k \), the functional
\[
J(q) = \int_{\mathbb{R}^k} \frac{g^2(\mathbf{x})}{q(\mathbf{x})}\, d\mathbf{x}, 
\qquad q \in \mathcal{D}_k.
\]
is uniquely minimised by \( g \) a.e.}
\end{lemma}

\begin{proof}

We apply Cauchy--Schwarz inequality:
\[
J(g) = 1 = \left( \int g(\mathbf{x})\, d\mathbf{x} \right)^2
= \left( \int \frac{g(\mathbf{x})}{q(\mathbf{x})} \, q(\mathbf{x})\, d\mathbf{x} \right)^2
\le \int \left( \frac{g(\mathbf{x})}{q(\mathbf{x})} \right)^2 q(\mathbf{x})\, d\mathbf{x}
= J(q),
\]

Equality in the Cauchy--Schwarz  inequality holds if and only if the argument of the convex function 
is constant almost everywhere with respect to \(q(\mathbf{x})\,d\mathbf{x}\). 
That is, there exists a constant \(c>0\) such that \(g(\mathbf{x})/q(\mathbf{x}) = c\) a.e. 
Since both \(g\) and \(q\) integrate to one, necessarily \(c=1\), and hence \(q = g\) almost everywhere.
\end{proof}

\subsubsection*{First step: optimisation of \( q_{-u \mid u} \)}

We optimise the conditional component $q_{-u \mid u}$ of $q$ for a fixed $q_u$. 
Since $\sigma^2_{\mathrm{opt},u}(q)$ in Eq.~(12) depends on $q_{-u \mid u}$ 
only through the term
$
\mathbb{E}_{p}\!\left[w_{-u\mid u}(\mathbf{X}_{-u}\mid\mathbf{X}_u)\,\phi^2(\mathbf{X}) \mid \mathbf{X}_u\right],
$
which appears in the definition of $\nu_q(\mathbf{X}_u)$, minimising this conditional expectation is sufficient. Applying Lemma~\ref{lem:variational_minimiser} conditionally on $\mathbf{X}_u$ yields the optimal conditional density
\begin{equation}
q^*_{-u \mid u}(\mathbf{x}_{-u} \mid \mathbf{x}_u)
= \frac{p_{-u \mid u}(\mathbf{x}_{-u}\mid \mathbf{x}_u)\,\phi(\mathbf{x})}
       {\mathbb{E}_{p}[\phi(\mathbf{X}) \mid \mathbf{X}_u=\mathbf{x}_u]},
\label{eq:q_star_conditional}
\end{equation}
where the denominator ensures normalisation.

With this choice, the variance contribution when 
$q_{-u \mid u} = q^*_{-u \mid u}$ becomes 
$\big( \mathbb{E}_{p}[ \phi(\mathbf{X}) \mid \mathbf{X}_u ] \big)^2$, 
instead of 
$\mathbb{E}_{p}[ \phi^2(\mathbf{X}) \mid \mathbf{X}_u ]$ 
when $\mathbf{X}_{-u} \mid \mathbf{X}_u \sim p_{-u \mid u}$ 
(i.e., $w_{-u \mid u} = 1$). 
A full derivation is provided in the Appendix.\\
Since the function \( x \mapsto x^2 \) is convex, Jensen’s inequality implies
\[
\big( \mathbb{E}_{p}[\phi(\mathbf{X}) \mid \mathbf{X}_u] \big)^2
\le \mathbb{E}_{p}[\phi^2(\mathbf{X}) \mid \mathbf{X}_u],
\]
showing that \(q^*_{-u\mid u}\) indeed reduces this part of the variance.

\subsubsection*{Second step: optimisation of the marginal \( q_u \)}

We now fix \( q_{-u \mid u} \), and seek to minimise the asymptotic variance over \( q_u \). 
For any such fixed conditional, the sampling density is of the form \( q = q_u \, q_{-u \mid u} \), 
and the corresponding variance reads
\[
\sigma_{\mathrm{opt},u}^2(q)
= \mathbb{E}_{p} \!\left[ w_u(\mathbf{X}_u)\, S(\mathbf{X}_u) \right] - \eta_u^2,
\]
where the integrand \( S(\mathbf{x}_u) \) depends on the choice of \( q_{-u \mid u} \). 
The optimal marginal minimising this expression is given by
\begin{equation}
q_u^*(\mathbf{x}_u)
= 
\frac{p_u(\mathbf{x}_u)\, \sqrt{S_{q_{-u \mid u}}(\mathbf{x}_u)}}
     {\mathbb{E}_{p}\!\left[\sqrt{S_{q_{-u \mid u}}(\mathbf{X}_u)}\right]},
\qquad \text{and} \qquad 
\sigma^2_{\mathrm{opt},u}(q_u^*)
=
\left( \mathbb{E}_{p}\!\left[\sqrt{S_{q_{-u \mid u}}(\mathbf{X}_u)}\right] \right)^2 - \eta_u^2.
\end{equation}

We distinguish below two particular cases:
\begin{itemize}
\item \textbf{Case A – Reference conditional:} \( q_{-u \mid u} = p_{-u \mid u} \)
\[
S_A(\mathbf{x}_u)
=
4\, m^2(\mathbf{x}_u)\,\nu_p(\mathbf{x}_u)
+m^4(\mathbf{x}_u)
\]

\item \textbf{Case B – Optimal conditional:} \( q_{-u \mid u} = q^*_{-u \mid u} \)
\[
S_B(\mathbf{x}_u)
=
4\, m^2(\mathbf{x}_u)\,\nu_{q^*}(\mathbf{x}_u)
+m^4(\mathbf{x}_u)
=
4\, m^2(\mathbf{x}_u)
\big(\mathbb{E}_{p}\!\left[\phi(\mathbf{X}) \mid \mathbf{X}_u = \mathbf{x}_u\right]\big)^2
-3m^4(\mathbf{x}_u),
\]
\end{itemize}

Both $S_A(\mathbf{x}_u)$ and $S_B(\mathbf{x}_u)$ are non-negative since they are built from conditional variances $\nu_p(\mathbf{X}_u)$ and $\nu_q(\mathbf{X}_u)$.

This satisfies the pointwise inequality 
\( S_B(\mathbf{x}_u) \le S_A(\mathbf{x}_u) \), as shown in the first step. 
The full chain of inequalities then reads
\[
\sigma^2_{\mathrm{opt},u}(q^*)
\le
\sigma^2_{\mathrm{opt},u}(q_u^* \, p_{-u \mid u})
\le
\sigma^2_{\mathrm{opt},u}(p),
\]
which follows from Jensen’s inequality:
\[
\left( \mathbb{E}_{p}\!\left[\sqrt{S_B(\mathbf{X}_u)}\right] \right)^2
\le 
\left( \mathbb{E}_{p}\!\left[\sqrt{S_A(\mathbf{X}_u)}\right] \right)^2
\le 
\mathbb{E}_{p}\!\left[S_A(\mathbf{X}_u)\right].
\]

Joint optimisation in Case B therefore achieves the greatest variance reduction, 
demonstrating the improved efficiency obtained when both components of \( q \) are optimised together.

\subsection{zero-variance optimal distribution}

In the idealised deterministic case, where all inputs are controllable and no external noise variable \(\mathbf{W}\) is present (\(Y=f(\mathbf{X})\)), it is theoretically possible to achieve zero variance in estimating \( \eta_u \).
This case is a direct simplification of the general framework, with the same two-step procedure but simpler expressions.
We still proceed in two steps, but the expressions simplify significantly. 
In particular, the asymptotic variance now takes the form
\begin{equation}
\sigma^2_{\mathrm{opt}, u}(q)
= 
\mathbb{E}_{p}\!\left[  
w_u(\mathbf{X}_u) 
\left(
4\, m^2(\mathbf{X}_u)\, 
\mathbb{E}_{p}\!\left[ 
w_{-u \mid u}(\mathbf{X}_{-u} \mid \mathbf{X}_u)\, Y^2 
\,\big|\, \mathbf{X}_u 
\right]
- 3\, m^4(\mathbf{X}_u) 
\right) 
\right] 
- \eta_u^2.    
\end{equation}\\

Applying Lemma~\ref{lem:variational_minimiser} sequentially yields
\begin{equation}
q^*_{-u \mid u}(\mathbf{x}_{-u} \mid \mathbf{x}_u)
= 
\frac{p_{-u \mid u}(\mathbf{x}_{-u} \mid \mathbf{x}_u)\, |f(\mathbf{x})|}
     {\mathbb{E}_{p}\!\left[ |f(\mathbf{X})| \,\big|\, \mathbf{X}_u = \mathbf{x}_u \right]},
\qquad
q_u^*(\mathbf{x}_u)
= 
\frac{p_u(\mathbf{x}_u)\, \sqrt{S_C(\mathbf{x}_u)}}
     {\mathbb{E}_{p}\!\left[\sqrt{S_C(\mathbf{X}_u)}\right]},
\label{eq:qstars-general}
\end{equation}
where
$$
S_C(\mathbf{x}_u)
= 
4\, m^2(\mathbf{x}_u)\,
\Big(\mathbb{E}_{p}\!\left[ |f(\mathbf{X})| \,\big|\, \mathbf{X}_u = \mathbf{x}_u \right]\Big)^2
- 3\, m^4(\mathbf{x}_u).
\label{eq:SC-def}\\
$$

Hence, the optimal joint sampling distribution is 
$
q^*(\mathbf{x}) 
= 
q_u^*(\mathbf{x}_u)\, q_{-u\mid u}^*(\mathbf{x}_{-u}\mid \mathbf{x}_u),
$ from which we obtain   
\[
\sigma^2_{\mathrm{opt},u}(q^*)
\le
\sigma^2_{\mathrm{opt},u}(p).
\] 
For instance, this distribution does not yield a zero-variance estimator, since a direct plug-in calculation shows 
\[
\sigma^2_{\mathrm{opt}, u}(q^*)
= 
\left(\mathbb{E}_{p}\!\left[\sqrt{S_C(\mathbf{X}_u)}\right]\right)^2 - \eta_u^2
> 0.
\]

\noindent\textbf{Remark.}
When the model \( f(\mathbf{x}) \) takes both positive and negative values, 
the associated optimal variance is strictly positive, so the zero-variance property cannot be achieved. 
If one nevertheless wishes to construct a zero-variance estimator, a classical alternative consists of separating 
the positive and negative contributions of \( f \), so that each part can be sampled using its own importance distribution. 
This idea, referred to as positivisation \cite{Owen2013}, 
relies on defining two auxiliary densities associated with the positive and negative parts of \( f \).

We decompose \( f \) into its positive and negative parts:
\[
f_+(\mathbf{x}) = \max(f(\mathbf{x}), 0),
\qquad
f_-(\mathbf{x}) = \max(-f(\mathbf{x}), 0),
\qquad
f(\mathbf{x}) = f_+(\mathbf{x}) - f_-(\mathbf{x}).
\]

Two auxiliary densities \( q_+ \) and \( q_- \) are then defined to be positive wherever 
\( p f_+ \) and \( p f_- \) are positive, respectively. 
Independent samples \( \mathbf{X}_{i,+} \sim q_+ \) and \( \mathbf{X}_{i,-} \sim q_- \) 
are then combined with opposite signs, yielding an unbiased estimator. 
In the ideal case where each component is sampled according to its fully adapted optimal density, the resulting estimator achieves zero variance.

If \( f(\mathbf{x}) \ge 0 \) for all \( \mathbf{x} \), then 
\(\mathbb{E}_{p}[|f(\mathbf{X})| \mid \mathbf{X}_u] = m(\mathbf{X}_u)\) 
and \( S_C(\mathbf{x}_u) = m^4(\mathbf{x}_u) \), so that
\begin{equation}
q^*_{-u \mid u}(\mathbf{x}_{-u} \mid \mathbf{x}_u)
= 
\frac{p_{-u \mid u}(\mathbf{x}_{-u} \mid \mathbf{x}_u)\, f(\mathbf{x})}
     {m(\mathbf{x}_u)},
\qquad
q_u^*(\mathbf{x}_u)
= 
\frac{p_u(\mathbf{x}_u)\, m^2(\mathbf{x}_u)}{\eta_u}.
\end{equation}
The resulting fully adapted sampling density is
\begin{equation}
q_u^*(\mathbf{x}) 
= 
\frac{p(\mathbf{x})\, f(\mathbf{x})\, m(\mathbf{x}_u)}{\eta_u},
\qquad \text{and} \qquad
\sigma^2_{\mathrm{opt}, u}(q^*) = 0.
\end{equation}
In theory, the density \( q^*(\mathbf{x}) \) yields a zero-variance estimator,  
but this ideal case cannot be achieved in practice 
since the optimal distribution depends on the unknown target quantity \( \eta_u \),
which brings us back to the classical importance sampling problem of approximating it.  
Adaptive procedures can be employed, including nonparametric methods \cite{zhang1996}  
and parametric ones such as the cross-entropy method \cite{papaioannou2019improved, deboer2005tutorial}.

\section{Numerical Experiments}
\label{sec:experiments}

\subsection{Analytical illustration of variance reduction}
In this section, we illustrate the effectiveness of an importance sampling transformation in reducing the variance of $\eta_u$ estimators. As a benchmark example, we consider the classical Sobol’ g-function \cite{Saltelli2000}:
\begin{equation}
    g(X_1, \dots, X_k) = \prod_{i=1}^k \frac{|4X_i - 2| + a_i}{1 + a_i},
\end{equation}
where $X_i \sim \mathcal{U}[0,1]$ are independent inputs. 
This function is a standard benchmark in sensitivity analysis due to its nonlinear, non-monotonic structure and the availability of analytical Sobol’ indices. 
We consider the 3-dimensional case $(k = 3)$ with $a_i = i$ and focus on the second-order index $S_{\{1,2\}}$, or equivalently $\eta_{\{1,2\}}$.For this part, one could in principle rely on kernel estimators \cite{DoksumSamarov1995,DaVeiga2013,DaVeiga2024}, but since all quantities are analytically tractable, we rely on their closed forms and focus on the impact of optimal importance densities on variance reduction. The third input $X_3$ is treated as an uncontrolled variable $W$. We consider Case A in Section 3.3, where $-u = \emptyset$. In this setting, all quantities involved in the expression of the optimal asymptotic variance can be computed analytically. Specifically, the optimal density $q^*$ to estimate $\eta_{\{1,2\}}$ is given by
\[
q^*(x_1, x_2) = \frac{p(x_1, x_2) \cdot \sqrt{S(x_1, x_2)}}{\mathbb{E}_{p}[\sqrt{S(X_1, X_2)}]},
\quad \text{with} \quad
S(x_1,x_2) = 4\, m^2(x_1,x_2)\, \phi^2(x_1,x_2)  - 3\, m^4(x_1,x_2).
\]
The conditional expectation $m(x_1,x_2) = \mathbb{E}[g(X_1, X_2, X_3) \mid X_1 = x_1, X_2 = x_2]$, and the squared term $\phi^2(x_1,x_2) =  \mathbb{E}[g^2(X_1,X_2,X_3)\mid X_1 = x_1, X_2 = x_2]$, can be derived analytically as $$ m(x_1,x_2) = \frac{(|4x_1 - 2| + 1)(|4x_2 - 2| + 2)}{6}, \quad
\phi^2(x_1,x_2) = \frac{99}{96} \cdot m^2(x_1,x_2).$$
This leads to the simplified expression
$$ S(x_1,x_2) = \frac{9}{8} \cdot m^4(x_1,x_2),
\quad \text{so that} \quad
q^*(x_1, x_2) = \frac{ m^2(x_1,x_2)}{\eta_{\{1,2\}}}.$$

\noindent\textbf{Remark.}
Although our example involves an uncontrolled variable \( W = X_3 \),
its particular structure satisfies 
\(\phi^2(x_u) \propto m^2(x_u)\),
so that the function \(S_A(x_u)\) simplifies to 
\(S_A(x_u) = c\, m^4(x_u)\) for some constant \(c>0\).
As a result,  \(q_u^*(x_u)\) retains the same analytical form
as in the noise-free case, namely \( q_u^*(x_u) \propto m^2(x_u) \).
This coincidence is specific to the g-function and does not hold in general models.

We first plot the shape of $q^*$ over \([0,1]^2\),  and then track how the coefficient of variation evolves when gradually moving from \( p \) to the optimal one $q^*$ via the interpolated family
$$
q_t(x_1, x_2) = (1 - t)\, p(x_1, x_2) + t\, q^*(x_1, x_2), \quad t \in [0,1].
$$
The coefficient of variation shown in Figure 1 is derived from the asymptotic variance under IS. It is defined as
$$
\mathrm{CV}_{\mathrm{opt}}(t)
= \frac{\sigma_{\mathrm{opt},u}(q_t)}{\eta_{\{1,2\}}},
\quad \text{with} \quad
\sigma^2_{\mathrm{opt},u}(q_t)
= \mathbb{E}_{p} \left[ \frac{S(X_1, X_2)}{q_t(X_1, X_2)} \right]
- \eta_{\{1,2\}}^2.
$$
Since all quantities are available analytically, $\mathrm{CV}_{\mathrm{opt}}(t)$ can be estimated using simple Monte Carlo, by sampling $X \sim p$ and evaluating the integrand directly.

\begin{figure}[H]
    \centering
    \begin{subfigure}[t]{0.43\textwidth}
        \centering
        \includegraphics[width=\linewidth]{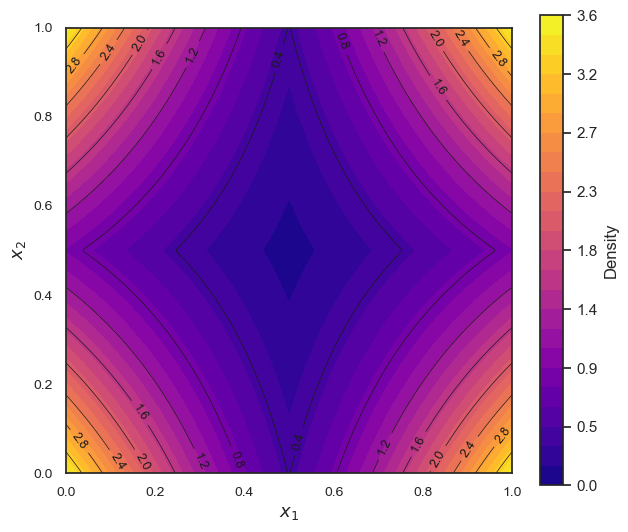}
        \caption{Optimal density \( q^*(x_1, x_2) \)}
    \end{subfigure}
    \hfill
    \begin{subfigure}[t]{0.53\textwidth}
        \centering
        \includegraphics[width=\linewidth]{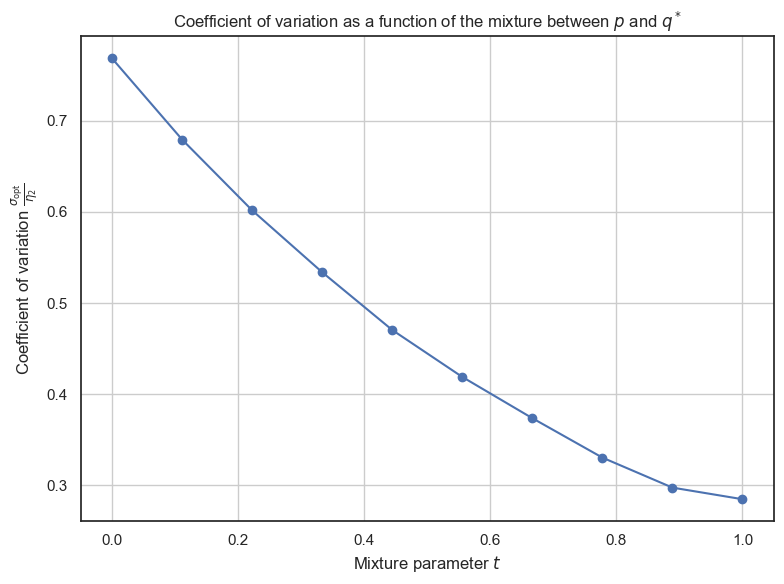}
        \caption{Coefficient of variation as a function of \( t \)}
    \end{subfigure}
    \caption{Left: optimal sampling density $q^*$. Right: relative dispersion of the estimator decreases as we move from \( p \) to the optimal \( q^* \).}
    \label{fig:q-optimality}
\end{figure}

\subsubsection*{Parametric optimisation}
An alternative strategy consists of introducing a parametric family of auxiliary sampling densities for $(X_1, X_2)$, defined by
$$
q_{\alpha,\beta}(x_1,x_2) = \mathrm{Beta}(x_1;\alpha_1,\beta_1) \cdot \mathrm{Beta}(x_2;\alpha_2,\beta_2),
$$
where the shape parameters $\alpha_i, \beta_i$ are positive real numbers, i.e., $\alpha_i, \beta_i \in \mathbb{R}_+^*$ for $i=1,2$. We restrict ourselves to the symmetric case $\alpha_1 = \alpha_2 = \alpha$ and $\beta_1 = \beta_2 = \beta$, so that the reference distribution corresponds to the uniform one at $(\alpha,\beta) = (1,1)$. The goal is to determine the optimal parameter pair $(\alpha^*, \beta^*) \in (\mathbb{R}_+^*)^2$ that minimises the asymptotic variance, namely
$$
(\alpha^*, \beta^*) = \arg\min_{(\alpha,\beta) \in (\mathbb{R}_+^*)^2} \sigma^2_{\mathrm{opt},u}(q_{\alpha,\beta}).
$$
The following plot shows how the asymptotic variance changes across the Beta family and highlights the optimal region.

\begin{figure}[H]
    \centering
    \includegraphics[width=0.6\textwidth]{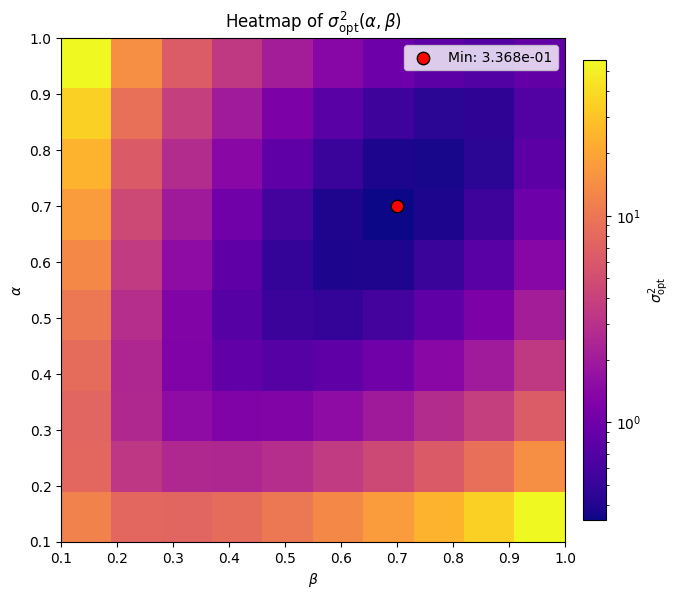}
    \caption{Asymptotic variance \(\sigma_{\mathrm{opt}}^2(\alpha, \beta)\) over the Beta parameter grid (symmetric case) displayed on a logarithmic scale.}
    \label{fig:heatmap_sigma_opt_log}
\end{figure}

The optimal configuration is found numerically around \((\alpha, \beta) \approx (0.7,\, 0.7)\). This yields a minimal variance of approximately \(0.3368\), corresponding to a reduction of about \(60\%\) compared to the uniform case \((\alpha = \beta = 1)\). The optimal parameters are marked in red on the heatmap.

\subsection{Sampling and estimation example}
To move closer to a realistic setting, we now consider an estimation scenario where the target quantity \( \eta_{\{1,2\}} \) is unknown and must be estimated from the observed data. We work with the same Sobol’ g-function and replace analytical quantities by Monte Carlo estimates while using the oracle \( m(\mathbf{x}) \) to isolate the impact of $q^*$.

We compare two strategies:
\begin{itemize}
    \item \textbf{Reference case:} samples drawn from the uniform input \( p \). In this case, we estimate
    \[
    \widehat{\eta}_{\{1,2\}}
    =
    \frac{1}{n} \sum_{i=1}^n (2Y_i - m(\mathbf{x}_i))\, m(\mathbf{x}_i),
    \qquad \mathbf{X}_i \sim p.
    \]

    \item \textbf{Importance sampling case:} samples are drawn from an empirical estimate of the optimal density $q^* \propto m^2(x_1,x_2)\,\mathbf{1}_{[0,1]^2}(x)$, as shown in the first part, denoted by \( \widehat{q}_\theta \).
    The corresponding estimator reads
    \[
    \widehat{\eta}_{\{1,2\}}
    =
    \frac{1}{n} \sum_{i=1}^n (2Z_i - m_Z(\mathbf{x}_i))\, m_Z(\mathbf{x}_i),
    \qquad
    Z_i=\frac{Y_i}{\sqrt{\widehat{q}_\theta(\mathbf{x}_i)}},\quad
    m_Z(\mathbf{x}_i)=\frac{m(\mathbf{x}_i)}{\sqrt{\widehat{q}_\theta(\mathbf{x}_i)}},\quad
      \mathbf{X}_i \sim \widehat{q}_\theta .
    \]
\end{itemize}

The Beta parametric family identified numerically in the first example can still be used as a simple baseline. In the following, however, we focus on a more flexible parametric approximation of \( q^* \). Since \( q^*(x_1,x_2) \propto m^2(x_1,x_2) \), its evaluation requires calls to the model
\( g(x_1,x_2,X_3) \) to compute the conditional mean \(m(x_1,x_2)\); consequently, the generation of samples approximately following \( q^* \) depends on model evaluations.\\

A preliminary sample \( (\mathbf{x}^{(i)})_{i=1}^{N_{\mathrm{fit}}} \) is obtained by a simple acceptance–rejection algorithm \cite{Forsythe1972}:
\[
\mathbf{x}^{(i)} \sim \mathrm{Unif}([0,1]^2), \qquad
u^{(i)} \sim \mathrm{Unif}(0,M), \qquad
\text{accept if } u^{(i)} \le m^2(\mathbf{x}^{(i)}),
\]
with \( M = \sup_{[0,1]^2} m^2(\mathbf{x})\). The accepted points follow \( q^* \) up to normalisation and are used to fit \(q_\theta\).

\begin{figure}[H]
    \centering
    \includegraphics[width=0.6\textwidth]{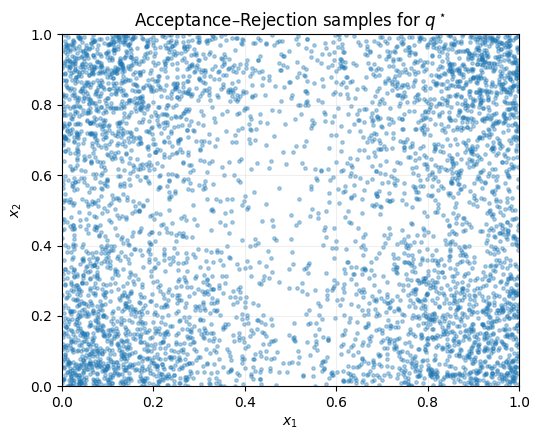}
    \caption{Samples generated via rejection sampling following \( q^* \).}
    \label{fig:qstar_samples}
\end{figure}

To preserve the support \([0,1]^2\) while retaining flexibility, we use a \emph{logit–Gaussian mixture model} (logit–GMM), where 
\(\mathrm{logit}(x) = \log\!\tfrac{x}{1-x}\) and \(\mathrm{logistic}(z) = (1 + e^{-z})^{-1}\).
Each pilot point \(\mathbf{x} = (x_1, x_2)\) is mapped to the unbounded space
\(\mathbf{z} = (z_1, z_2)\) through \( z_j = \mathrm{logit}(x_j) \), and a Gaussian mixture is then fitted:
\[
q_Z(\mathbf{z}; \theta)
= \sum_{k=1}^{K} \pi_k\, \mathcal{N}(\mathbf{z}; \mu_k, \Sigma_k),
\qquad
\pi_k > 0, \quad \sum_{k=1}^{K} \pi_k = 1,
\]
with parameters \( \theta = \{\pi_k, \mu_k, \Sigma_k\}_{k=1}^{K} \).

By the change-of-variables formula, this induces a properly normalised expression on \([0,1]^2\):
\[
q_\theta(\mathbf{x})
= q_Z\!\bigl(\mathrm{logit}(\mathbf{x}); \theta\bigr)
  \prod_{j=1}^2 \frac{1}{x_j(1-x_j)},
  \qquad \mathbf{x} \in (0,1)^2.
\]

\begin{figure}[H]
    \centering
    \includegraphics[width=0.6\textwidth]{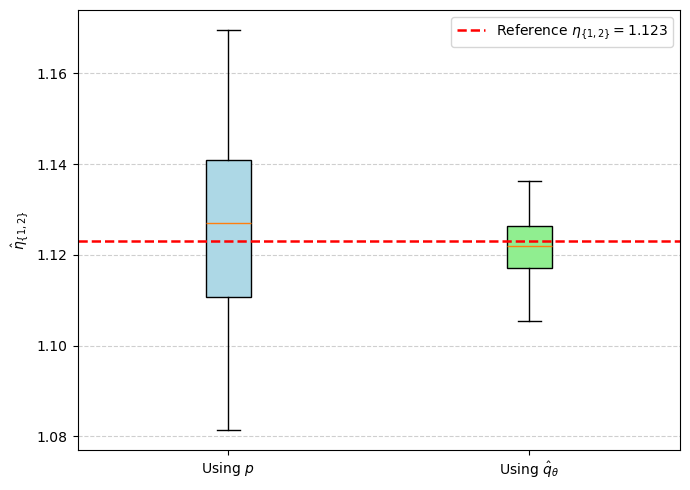}
    \caption{Dispersion of \( \widehat{\eta}_{\{1,2\}} \) under \( p \) and the learned \( \widehat{q}_\theta \).}
    \label{fig:boxplot_eta}
\end{figure}

The fitted \( \widehat{q}_\theta \) maximises the likelihood of the pilot sample under the GMM using the EM algorithm. 
This formulation ensures proper normalisation over the unit square and avoids boundary bias. 
Sampling under $\widehat{q}_\theta$ is straightforward: we draw 
$\mathbf{Z} \sim q_Z(\cdot;\widehat{\theta})$, set 
$\mathbf{X}=\mathrm{logistic}(\mathbf{Z})$, and use the same transformed density 
$\widehat{q}_\theta(\mathbf{x})$ to compute importance weights. 
These samples $\mathbf{X}_i \sim \widehat{q}_\theta$ are then used to evaluate the IS estimator. The number of mixture components is fixed to $K=4$, which provides sufficient flexibility to approximate $q^*$ while remaining numerically stable for $N_{\mathrm{fit}}=1000$.
For fairness, both approaches are allocated the same total number of model evaluations. 
The density $\widehat{q}_{\theta}$ is learned from $N_{\mathrm{fit}} = 1000$ evaluations, 
and the final estimation uses an additional $n = 1000$ evaluations, for a total cost of $2000$. 
Therefore, the estimator under $p$ is computed with $n = 2000$ samples so that both methods 
rely on the same overall computational budget.
The experiment is repeated $100$ times, and the resulting empirical distributions 
are summarised using boxplots.
Figure~\ref{fig:boxplot_eta} compares the dispersion of the estimates, showing a clear variance reduction under the learned \( \widehat{q}_\theta \).

This experiment mimics a practical learning-based IS setting, where the optimal distribution $q^*$ is unknown and approximated from data, leading to noticeable variance reduction.

\subsection{Distributional Sensitivity via Reverse Importance Sampling}
We study the sensitivity of the quantity $ \eta_u(\theta) = \mathbb{E}_{p_\theta}\left[ m^2(\mathbf{X}_u) \right],$
with respect to changes in the input distribution \(p_\theta\), based on a single dataset \((\mathbf{X}^{(i)}, Y^{(i)})_{i=1}^n\) generated under \(p\). This falls into the framework of reverse importance sampling. 
The rank estimator \cite{Gamboa2022}, based on the sequential correlation of outputs, takes as input a sequence $Y^{(i)}$ sorted with respect to the associated inputs $\mathbf{X}_u^{(i)}$, and converges to $\eta_u$ via the average of consecutive products. 
By using the results from Section~3.1, the same estimator can be extended to compute $\eta_u(\theta)$ under any distribution $p_\theta$. 
This formulation allows extending the rank-based estimator of \cite{Gamboa2022}
to any perturbed distribution \(p_\theta\) through reverse importance sampling and remains consistent.
Its asymptotic behavior follows directly from the consistency result established in Proposition~3.2 of \cite{Gamboa2022}.
We note that in this section we work with the rank estimator (and not the lagged–rank estimator), which is not asymptotically efficient in terms of variance. Since our purpose here is only to study the sensitivity of the Sobol’ indices under different input distributions, the optimal densities derived in Section~3 are not involved. For completeness, we nevertheless derive in Appendix~B the optimal sampling densities associated with the rank estimator itself.
For each \( \theta \), we define the reweighted output and the rank-based estimator becomes
\begin{equation}
\widehat{\eta}^{\mathrm{rank}}_u(\theta)
=
\frac{1}{n-1}\sum_{i=1}^{n-1}
Z_{\theta,u}^{(i)}\, Z_{\theta,u}^{(i+1)},
\qquad
Z_{\theta,u}
=
\sqrt{\frac{p_{\theta,u}(\mathbf{X}_u)}{p_u(\mathbf{X}_u)}}
\,
\frac{p_{\theta,- u}(\mathbf{X}_{- u} \mid \mathbf{X}_u)}
     {p_{- u}(\mathbf{X}_{- u} \mid \mathbf{X}_u)}
\,
Y.
\label{eq:rank-based-estimator}
\end{equation}

Hence, from a single dataset generated under \( p \), one can efficiently explore the sensitivity of \( \eta_u \) (and thus the associated Sobol' index \( S_u \)) across the entire family of distributions \( p_\theta \), without performing any additional evaluations of the black-box model.
\subsection*{Illustrative test case}
\paragraph{Physical model}  
We consider a black-box simulator based on the IGLOO2D code developed at ONERA \cite{Trontin2017}, which returns the ice-accretion thickness along an airfoil profile. Each run of the simulator takes as input a vector of four physical variables $\mathbf{X} = (X_1, X_2, X_3, X_4) \in \mathbb{R}^4$, initially sampled independently from uniform distributions defined over the intervals $\mathbf{x}_{\mathrm{low}} = (3.5,\;65,\;261.15,\;4 \times 10^{-4})$ and $\mathbf{x}_{\mathrm{up}} = (4.5,\;75,\;265.15,\;6 \times 10^{-4})$, which respectively correspond to $X_1$ (angle of attack, AoA), $X_2$ (freestream velocity, $U_\infty$), $X_3$ (ambient temperature, $T_\mathrm{inf}$), and $X_4$ (liquid water content, LWC). A design of experiments of size $N = 2500$ has been generated, yielding input samples $\mathbf{X}^{(i)} \in \mathbb{R}^4$ for $i = 1,\dots,N$, and corresponding outputs $Y^{(i)} \in \mathbb{R}$, which represent the accretion thickness at a given location along the airfoil. All input variables $X_j$ are first normalised to the unit interval via $X_j^{\mathrm{std}} := (X_j - x_{\mathrm{low},j}) / (x_{\mathrm{up},j} - x_{\mathrm{low},j})$. To explore the sensitivity of this output to changes in the input distributions, we consider alternative input distributions $p_{\theta_j}(X_j)$ built from rescaled $\operatorname{Beta}(\theta_j)$ marginals on the same support, where varying $\theta_j = (\alpha_j, \beta_j)$ allows smooth deviations from the uniform case $(1,1)$. For any configuration $\boldsymbol{\theta} = (\theta_1, \dots, \theta_4)$, the corresponding Sobol' index $S_u(\boldsymbol{\theta})$ or contrast $\eta_u(\boldsymbol{\theta})$ can be estimated from the same dataset by importance weighting, using the factor
\begin{equation}
w_{\boldsymbol{\theta}}(\mathbf{X}) = \prod_{j=1}^{4} \frac{\operatorname{Beta}_{\theta_j}(X_j^{\mathrm{std}})}{\operatorname{Unif}[0,1](X_j^{\mathrm{std}})}.  \end{equation}
With the rank estimator in Eq.~(20), the quantity $\eta_u(\boldsymbol{\theta}) = \mathbb{E}_{p_{\boldsymbol{\theta}}}[m_u^2(\mathbf{X}_u)]$ can be estimated for any configuration of Beta parameters $\boldsymbol{\theta} = (\theta_1, \dots, \theta_4)$. Two perturbation strategies are considered: (i) \emph{marginal perturbation}, where only the pair $\theta_j = (\alpha_j, \beta_j)$ associated with a single input $X_j$ is varied while the remaining marginals remain uniform; and (ii) \emph{global perturbation}, where several or all input marginals are varied simultaneously.
To keep the presentation simple, we focus on marginal perturbations. Specifically, we assess the distributional sensitivity of the contrast $\eta_u$ with $u = \{4\}$ by varying $\theta_j$ for one input variable at a time, keeping all others unchanged. This yields a family of perturbed estimates $\eta_4(\theta_j)$, each revealing how the index $\eta_4$ responds to changes in the distribution of a specific input $X_j$, for $j = 1, \dots, 4$. All computations are based on the same reference dataset.

\begin{figure}[H]
  \centering

  \begin{subfigure}[t]{0.44\textwidth}
    \centering
    \includegraphics[width=\linewidth]{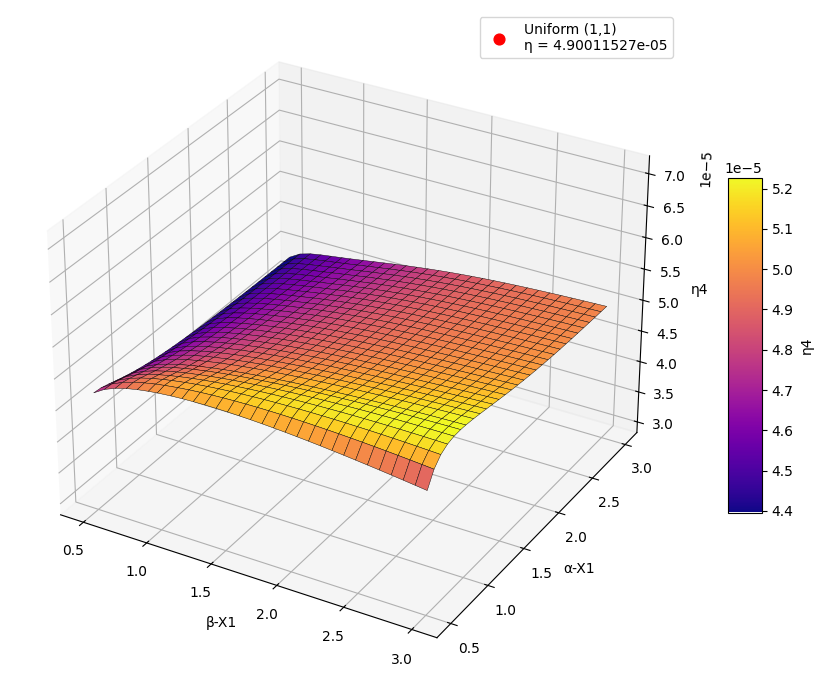}
    \caption{Varying \( \alpha_1, \beta_1 \) of \( X_1 \)}
  \end{subfigure}
  \hfill
  \begin{subfigure}[t]{0.44\textwidth}
    \centering
    \includegraphics[width=\linewidth]{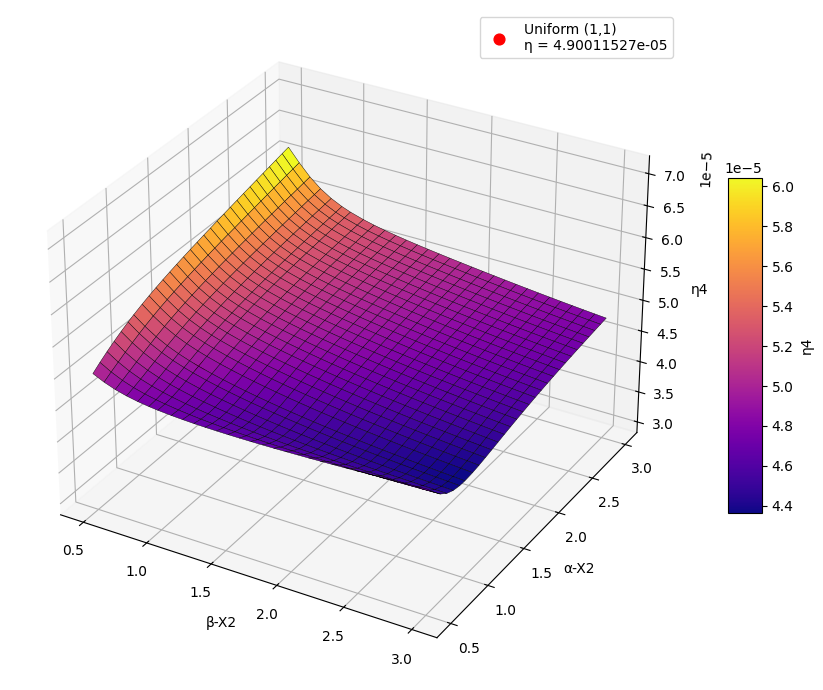}
    \caption{Varying \( \alpha_2, \beta_2 \) of \( X_2 \)}
  \end{subfigure}

  \medskip

  \begin{subfigure}[t]{0.44\textwidth}
    \centering
    \includegraphics[width=\linewidth]{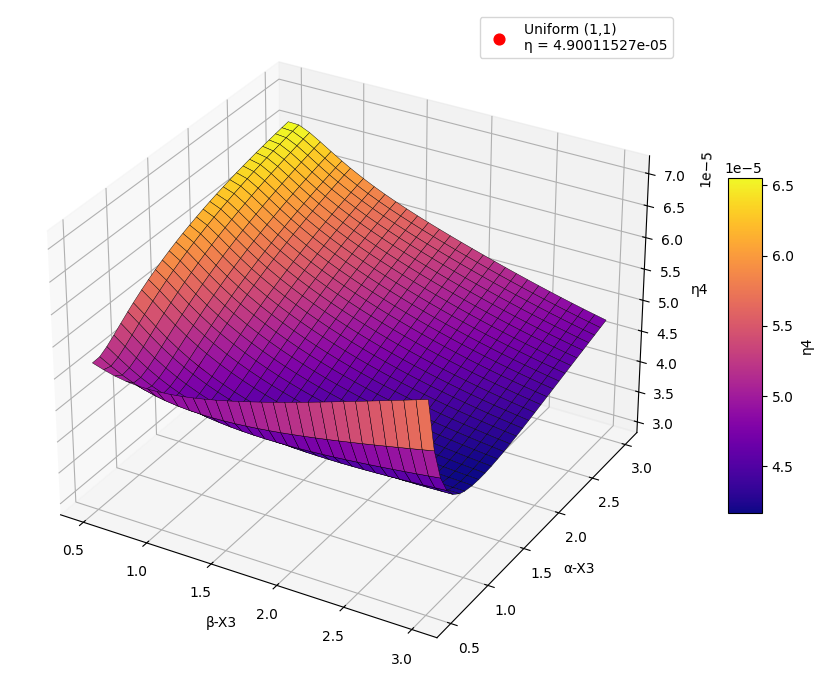}
    \caption{Varying \( \alpha_3, \beta_3 \) of \( X_3 \)}
  \end{subfigure}
  \hfill
  \begin{subfigure}[t]{0.44\textwidth}
    \centering
    \includegraphics[width=\linewidth]{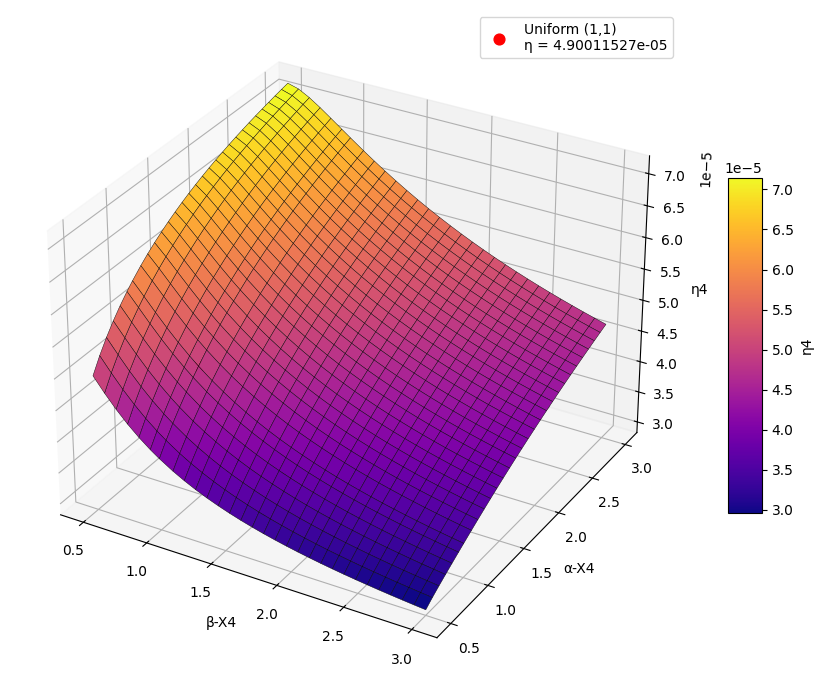}
    \caption{Varying \( \alpha_4, \beta_4 \) of \( X_4 \)}
  \end{subfigure}

  \caption{Effect of marginal distribution changes on \( \eta_4 \):  
           each subplot shows how \( \eta_4 \) evolves as the Beta parameters of one input \( X_j \) are varied, while others are kept uniform.  
           All curves are computed from the same baseline sample (\( N = 2500 \)).}
  \label{fig:eta-sensitivity}
\end{figure}

These results show that the sensitivity of \( \eta_4(\theta_j) \) strongly depends on the shape of the input distribution of \( X_j \), which varies through the two parameters \( \alpha_j \) and \( \beta_j \) of the Beta distribution, with different behaviours across the four inputs.
 The strongest variations are observed for the pair \( (\alpha_4, \beta_4) \), which is expected since \( \eta_4 \) quantifies the influence of $X_4$. Overall, these results confirm that the proposed given‐data IS framework offers a powerful and cost-effective means to perform detailed distributional local sensitivity analyses of Sobol’ indices.

\section{Conclusion}
\label{sec:conclusions}

In this paper, we develop a new framework for the estimation of Sobol' indices using importance sampling.  
We demonstrate that $\mathbb{E}\left[ m_u^2(\mathbf{X}_u) \right]$ can be estimated from any sampling distribution that is different from the reference one. In particular, we show that by choosing suitable auxiliary densities, the asymptotic variance of efficient estimators can be further reduced, as changing the sampling distribution modifies the efficiency bound itself. In some cases, zero-variance estimation can theoretically be achieved; however, this is not feasible in practice since the optimal densities depend on unknown quantities. In realistic settings where model evaluations are available, one may nevertheless approximate these optimal densities, which naturally leads to significant variance reduction.
We also show that combining our method with the rank-based estimator \cite{Gamboa2022} is particularly useful for studying the distributional sensitivity of Sobol' indices, especially in the given-data setting where samples can be reused efficiently through the reverse importance sampling mechanism \cite{Lemaitre2015}. The first limitation of our present work is that the optimal densities depend on the code, so additional computational effort is required. 
Another limitation is that we have proved the existence of one optimal distribution for every Sobol' index, so it may be computationally expensive to estimate all these optimal densities, which increases the numerical effort. 
In a related direction, \cite{DemangeChryst2023} investigate how several expectations under different target distributions can be estimated from a single Monte Carlo sample, using an adaptive importance sampling scheme with shared control variates. This suggests that multi–target designs could also be explored to construct sampling densities efficient for several Sobol’ quantities simultaneously.
Future work may focus on constructing joint optimal densities for the simultaneous estimation of all $\eta_u$, rather than treating them independently.

\section*{Acknowledgement}

We gratefully acknowledge our ONERA colleagues — Lokman Bennani, Ghislain Blanchard, Maxime Bouyges, and Patricia Klotz — for providing the ice-accretion dataset. We also thank the reviewers for their valuable remarks during the revision process.

\section*{Supplementary Materials}
All simulation codes for the academic experiments (Section~4)  
are available in the Supplementary Materials and on GitHub:  
\[
\text{{https://github.com/HaythemBC/Importance-sampling-for-Sobol-indices-estimation}}.
\]

\section*{Data Availability Statement}
All data used in the academic case studies are generated by simulation  
and can be fully reproduced using the associated code.  
The real-world dataset of Section~4.3 is confidential and cannot be shared publicly,  
but can be provided to reviewers upon reasonable request, in accordance with the data-sharing agreement of the contributing institution.

\appendix
\section{Lemma 3.1 and other derivations} 

\textbf{Proof of Lemma 3.1.}
We compute the conditional expectation of \( Z_u \) given \( X_u = x_u \), and deduce \( \eta_u \) by taking successive expectations under \( q \).

\[
\begin{aligned}
\mathbb{E}_{q}[Z_u \mid \mathbf{X}_u]
&= \sqrt{w_u(\mathbf{X}_u)} \cdot \mathbb{E}_{q} \left[ w_{-u}(\mathbf{X}_{-u} \mid \mathbf{X}_u) \cdot \mathbb{E}[Y \mid \mathbf{X}] \,\middle|\, \mathbf{X}_u \right] \\[6pt]
&= \sqrt{w_u(\mathbf{X}_u)} \cdot \mathbb{E}_{p} \left[ \mathbb{E}[Y \mid \mathbf{X}] \,\middle|\, \mathbf{X}_u \right] \\[6pt]
&= \sqrt{w_u(\mathbf{X}_u)} \cdot m(\mathbf{X}_u), 
\end{aligned}
\]
and
\[
\mathbb{E}_{q}[\mathbb{E}_{q}[Z_u \mid \mathbf{X}_u]^2 ] = \mathbb{E}_{q}[w_u(\mathbf X_u) m^2(\mathbf{X}_u) ] = \mathbb{E}_{p}[m^2(\mathbf{X}_u) ] 
\]
This completes the proof.

\textbf{Derivation of the optimal variance under \(p\) (Eq.~(12)).}
From Lemma~3.1, we recall
\[
\mathbb{E}_q[Z_u \mid X_u] = \sqrt{w_u(X_u)}\, m(X_u),
\qquad
\bigl(\mathbb{E}_q[Z_u\mid X_u]\bigr)^2 = w_u(X_u)m^2(X_u).
\]
\[
\nu_q(X_u)
 := \mathbb{E}_p[w_{-u\mid u}\,\phi^2(X)\mid X_u] - m^2(X_u).
\]\\
We now compute the conditional second moment. Using the definition of
\(Z_u\) and
\(\phi^2(x)\), we obtain
\[
\begin{aligned}
\mathbb{E}_q[Z_u^2\mid X_u]
 &= \mathbb{E}_q\!\left[w_u(X_u) w_{-u\mid u}(X_{-u}\mid X_u) Y^2 
    \mid X_u\right] \\[2pt]
 &= w_u(X_u)\,
    \mathbb{E}_p\!\left[w_{-u\mid u}(X_{-u}\mid X_u)\,\phi^2(X)
    \mid X_u\right].
\end{aligned}
\]

Hence
\[
\operatorname{Var}_q(Z_u\mid X_u)
 = w_u(X_u)\Big(
     \mathbb{E}_p[w_{-u\mid u}\,\phi^2(X)\mid X_u]
     - m^2(X_u)
   \Big)
 = w_u(X_u)\,\nu_q(X_u),\]\\
Substituting all expressions 
into Proposition~3.2 yields
\[
\begin{aligned}
\sigma^2_{\mathrm{opt},u}(q)
 &= \mathbb{E}_q\!\left[
      4\,w_u(X_u)m^2(X_u)\, w_u(X_u)\nu_q(X_u)
      + w_u^2(X_u)m^4(X_u)
    \right]
    - \eta_u^2 \\[4pt]
 &= \mathbb{E}_q\!\left[
      w_u^2(X_u)\bigl(4 m^2(X_u)\nu_q(X_u) + m^4(X_u)\bigr)
    \right] - \eta_u^2.
\end{aligned}
\]

Using the usual identity
\(\mathbb{E}_q[w_u^2 g] = \mathbb{E}_p[w_u g]\), we finally obtain
\[
\sigma^2_{\mathrm{opt},u}(q)
 = \mathbb{E}_p\!\left[
     w_u(X_u)\bigl(4 m^2(X_u)\nu_q(X_u) + m^4(X_u)\bigr)
   \right] - \eta_u^2,
\]
which is Eq.~(12) of the main text.

\textbf{Computing optimal variance bounds.} 
Plugging $ q^*_{-u \mid u}(\mathbf{x}_{-u} \mid \mathbf{x}_u) $, we compute the quantity 
that contributes to the variance 
$\sigma^2_{\mathrm{opt}}(q = (p_u, q^*_{-u \mid u}))$. 
We obtain:

\begin{align*}
\mathbb{E}_{p} 
\left[ w(\mathbf{X}) \cdot \phi^2(\mathbf{X}) \mid \mathbf{X}_u \right]
&= \mathbb{E}_{p} 
\left[ \frac{p(\mathbf{X})}{q^*(\mathbf{X} \mid \mathbf{X}_u)} \cdot \phi^2(\mathbf{X}) \,\middle|\, \mathbf{X}_u \right] \\[0.6em]
&= \mathbb{E}_{p} 
\left[ \frac{\mathbb{E}_{p}[\phi(\mathbf{X}) \mid \mathbf{X}_u]}{\phi(\mathbf{X})} \cdot \phi^2(\mathbf{X}) \,\middle|\, \mathbf{X}_u \right] \\[0.6em]
&= \mathbb{E}_{p} 
\left[ \phi(\mathbf{X}) \cdot \mathbb{E}_{p}[\phi(\mathbf{X}) \mid \mathbf{X}_u] \,\middle|\, \mathbf{X}_u \right] \\
&= \left( \mathbb{E}_{p}[\phi(\mathbf{X}) \mid \mathbf{X}_u] \right)^2.
\end{align*}

Similarly for \(q^*(\mathbf{x}_u)\):
\[
\begin{aligned}
\mathbb{E}_{p} \left[ w(\mathbf{X}_u)\, S(\mathbf{X}_u) \right] - \eta_u^2
&= \mathbb{E}_{p} \left[ \frac{p(\mathbf{X}_u)}{q^*(\mathbf{X}_u)}\, S(\mathbf{X}_u) \right] - \eta_u^2 \\[0.6em]
&= \mathbb{E}_{p} \left[ \frac{ \mathbb{E}_{p}[\sqrt{S(\mathbf{X}_u)}] }{ \sqrt{S(\mathbf{X}_u)} }\, S(\mathbf{X}_u) \right] - \eta_u^2 \\[0.6em]
&= \mathbb{E}_{p} \left[ \mathbb{E}_{p}[\sqrt{S(\mathbf{X}_u)}] \cdot \sqrt{S(\mathbf{X}_u)} \right] - \eta_u^2 \\[0.6em]
&= \left( \mathbb{E}_{p}[\sqrt{S(\mathbf{X}_u)}] \right)^2 - \eta_u^2.
\end{aligned}
\]

\textbf{Derivation of the zero-variance formula}
We start from the general expression:
\[
\sigma^2_{\mathrm{opt}, u}(q)
= \mathbb{E}_{p} \left[
    w_u(\mathbf{X}_u) \cdot \left(
        4\, m^2(\mathbf{X}_u) \cdot \mathbb{E}_{p}\left[ w_{-u}(\mathbf{X}_{-u} \mid \mathbf{X}_u)\, \phi^2(\mathbf{X}) \mid \mathbf{X}_u \right]
        - 3\, m^4(\mathbf{X}_u)
    \right)
\right]
- \eta_u^2.
\]

Now plugging both optimal densities simultaneously.

$$
q_u^*(\mathbf{x}_u) = \frac{p_u(\mathbf{x}_u)\, m^2(\mathbf{x}_u)}{\eta_u}  \qquad
q_{- u}^*(\mathbf{x}_{- u} \mid \mathbf{x}_u) =
\frac{p_{- u}(\mathbf{x}_{- u})\, f(\mathbf{x})}{m^2(\mathbf x_u)}
$$

\[
\begin{aligned}
\sigma^2_{\mathrm{opt}, u}(q)
&= \mathbb{E}_{p} \left[
    w(\mathbf{X}) \cdot \left(
        4\, m^2(\mathbf{X}) \cdot \mathbb{E}_{p} \left[
            w(\mathbf{X}) \cdot \phi^2(\mathbf{X}) \mid \mathbf{X} \right]
        - 3\, m^4(\mathbf{X})
    \right)
\right]
- \eta_u^2 \\[0.6em]
&= \mathbb{E}_{p} \left[
    \frac{\eta_u}{m^2(\mathbf{X})} \cdot \left(
        4\, m^2(\mathbf{X}) \cdot \mathbb{E}_{p} \left[
            \frac{m(\mathbf{X})}{f(\mathbf{X})} \cdot f^2(\mathbf{X}) \mid \mathbf{X}
        \right]
        - 3\, m^4(\mathbf{X})
    \right)
\right]
- \eta_u^2 \\[0.6em]
&= \mathbb{E}_{p} \left[
    \frac{\eta_u}{m^2(\mathbf{X})} \cdot \left(
        4\, m^3(\mathbf{X}) \cdot m(\mathbf{X})
        - 3\, m^4(\mathbf{X})
    \right)
\right]
- \eta_u^2 \\[0.6em]
&= \mathbb{E}_{p} \left[
    \frac{\eta_u}{m^2(\mathbf{X})} \cdot m^4(\mathbf{X})
\right]
- \eta_u^2 \\[0.6em]
&= \eta_u \cdot \mathbb{E}_{p} \left[ m^2(\mathbf{X}) \right]
- \eta_u^2\\[0.6em]
&= \eta_u^2 - \eta_u^2 = 0.
\end{aligned}
\]

\section*{Appendix B. Optimal sampling densities for the rank estimator}

In this appendix, we derive the optimal sampling densities for the rank-based estimator of 
$\eta_u = \mathbb{E}[m^2(\mathbf{X}_u)]$.  
Since this estimator is not asymptotically efficient, the optimal variance bound of Section~3 does not apply.  
Nevertheless, the same two-step optimisation over $q_{-u\mid u}$ and $q_u$ can be performed,
leading to sampling densities specifically adapted to its asymptotic variance.
Throughout, we use $m(\mathbf{x}_u)=\mathbb{E}[Y\mid \mathbf{X}_u=\mathbf{x}_u]$,
$\nu(\mathbf{x}_u)=\operatorname{Var}(Y\mid \mathbf{X}_u=\mathbf{x}_u)$, 
and $\phi^2(\mathbf{x})=\mathbb{E}[Y^2\mid \mathbf{X}=\mathbf{x}]$. We first recall the asymptotic variance under the reference distribution:
\[
n\,\operatorname{Var}(\widehat{\eta}_{\mathrm{rank},u})
=
4\,\mathbb{E}\big[m^2(\mathbf{X}_u)\,\nu(\mathbf{X}_u)\big]
+
\mathbb{E}\big[\nu^2(\mathbf{X}_u)\big]
+
\operatorname{Var}\big(m^2(\mathbf{X}_u)\big)
+
o(1),
\]
which can be rewritten as
\[
\sigma^2_{\mathrm{rank},u}(p)
=
\mathbb{E}_{p}\big[
4\,m^2(\mathbf{X}_u)\,\nu_p(\mathbf{X}_u)
+
m^4(\mathbf{X}_u)
+
\nu_p^2(\mathbf{X}_u)
\big]
-
\eta_u^2.
\]

Following exactly the derivation of Proposition~3.2, the variance under any density $q$ can be written as
\[
\sigma^2_{\mathrm{rank},u}(q)
=
\mathbb{E}_p\!\big[ w_u(\mathbf{X}_u)\,S_{\mathrm{rank},q}(\mathbf{X}_u)\big] - \eta_u^2,
\qquad \text{where} \qquad
S_{\mathrm{rank},q}(\mathbf{x}_u)
=
4\,m^2(\mathbf{x}_u)\,\nu_q(\mathbf{x}_u)
+
\nu_q^2(\mathbf{x}_u)
+
m^4(\mathbf{x}_u),
\]
and
\[
\nu_q(\mathbf{x}_u)
=
\mathbb{E}_p\!\left[
w_{-u\mid u}(\mathbf{X}_{-u}\mid \mathbf{X}_u)\,\phi^2(\mathbf{X})
\mid \mathbf{X}_u=\mathbf{x}_u
\right]
-
m^2(\mathbf{x}_u).
\]

As in Section~3, $S_{\mathrm{rank},q}(\mathbf{x}_u)$ depends on $q_{-u\mid u}$ only through $\nu_q(\mathbf{x}_u)$.
Hence minimising $S_{\mathrm{rank},q}(\mathbf{x}_u)$ reduces to minimising $\nu_q(\mathbf{x}_u)$,
and the optimal conditional density is exactly the same as in \eqref{eq:q_star_conditional}:
\[
q^*_{-u\mid u}(\mathbf{x}_{-u}\mid \mathbf{x}_u)
=
\frac{
p_{-u\mid u}(\mathbf{x}_{-u}\mid \mathbf{x}_u)\,\phi(\mathbf{x})
}{
\mathbb{E}_p[\phi(\mathbf{X})\mid \mathbf{X}_u=\mathbf{x}_u]
}.
\]
Under this density,
\[
\nu_{q^*}(\mathbf{x}_u)
=
\big(\mathbb{E}_p[\phi(\mathbf{X})\mid \mathbf{X}_u=\mathbf{x}_u]\big)^2 - m^2(\mathbf{x}_u).
\]

\paragraph{Case A and Case B}  
As in Section~3, we distinguish the two choices for the conditional component:
\[
S_A(\mathbf{x}_u)
=
4\,m^2(\mathbf{x}_u)\,\nu_p(\mathbf{x}_u)
+
\nu_p^2(\mathbf{x}_u)
+
m^4(\mathbf{x}_u),
\qquad
\text{when } q_{-u\mid u}=p_{-u\mid u},
\]
and
\[
S_B(\mathbf{x}_u)
=
4\,m^2(\mathbf{x}_u)\,\nu_{q^*}(\mathbf{x}_u)
+
\nu_{q^*}^2(\mathbf{x}_u)
+
m^4(\mathbf{x}_u),
\qquad
\text{when } q_{-u\mid u}=q^*_{-u\mid u}.
\]\\
Compared to Section~3, the only difference is the additional term $\nu^2$ in both $S_A$ and $S_B$. For a fixed conditional $q_{-u\mid u}$ (either Case~A or Case~B), the variance can be written as
\[
\sigma^2_{\mathrm{rank},u}(q)
=
\mathbb{E}_p\!\left[w_u(\mathbf{X}_u)\,S(\mathbf{X}_u)\right] - \eta_u^2,
\]
where $S$ denotes $S_A$ or $S_B$ depending on the chosen conditional component. As in Section~3, the optimal marginal is
\[
q_u^*(\mathbf{x}_u)
=
\frac{
p_u(\mathbf{x}_u)\,\sqrt{S(\mathbf{x}_u)}
}{
\mathbb{E}_p[\sqrt{S(\mathbf{X}_u)}]
}, \qquad \text{and} \qquad \sigma^2_{\mathrm{rank},u}(q_u^*, q_{-u\mid u})
=
\left(
\mathbb{E}_p[\sqrt{S(\mathbf{X}_u)}]
\right)^2
-\eta_u^2.
\]

\paragraph{Resulting optimal density.}
The optimal sampling density for the rank estimator is therefore
\[
q_u^*(\mathbf{x})
=
q_u^*(\mathbf{x}_u)\,
q^*_{-u\mid u}(\mathbf{x}_{-u}\mid \mathbf{x}_u),
\]
with the conditional component identical to Section~3, while the marginal depends on
$S_A$ or $S_B$ through the additional term $\nu^2$, which encodes the intrinsic sub-efficiency
of the rank estimator.

\paragraph{Deterministic case ($Y=f(\mathbf X)$).}
Assume that $f\ge 0$ and that there is no stochastic noise, i.e.\ $Y=f(\mathbf X)$.
Then $\phi^2(\mathbf x)=\mathbb E[Y^2\mid \mathbf X=\mathbf x]=f^2(\mathbf x)$ and thus $\phi(\mathbf x)=f(\mathbf x)$.
In particular,
\[
m(\mathbf x_u)=\mathbb E_p[f(\mathbf X)\mid \mathbf X_u=\mathbf x_u],
\qquad
\nu_p(\mathbf x_u)=\operatorname{Var}_p(f(\mathbf X)\mid \mathbf X_u=\mathbf x_u).
\]
The optimal conditional density becomes
\[
q^*_{-u\mid u}(\mathbf x_{-u}\mid \mathbf x_u)
=
\frac{p_{-u\mid u}(\mathbf x_{-u}\mid \mathbf x_u)\,f(\mathbf x)}
{\mathbb E_p[f(\mathbf X)\mid \mathbf X_u=\mathbf x_u]}
=
\frac{p_{-u\mid u}(\mathbf x_{-u}\mid \mathbf x_u)\,f(\mathbf x)}{m(\mathbf x_u)}.
\]
Moreover, under $q^*_{-u\mid u}$ we have
\[
\nu_{q^*}(\mathbf x_u)
=
\big(\mathbb E_p[f(\mathbf X)\mid \mathbf X_u=\mathbf x_u]\big)^2 - m^2(\mathbf x_u)
=0,
\]
so that $S_B(\mathbf x_u)=m^4(\mathbf x_u)$ and hence
\[
q_u^*(\mathbf x_u)=\frac{p_u(\mathbf x_u)\,m^2(\mathbf x_u)}{\eta_u}.
\]
Therefore, the resulting optimal sampling density is
\[
q_u^*(\mathbf x)
=
\frac{p(\mathbf x)\,f(\mathbf x)\,m(\mathbf x_u)}{\eta_u}.
\]


\end{document}